
\magnification1200
\input amstex.tex
\documentstyle{amsppt}

 \hsize=12.5cm
 \vsize=18cm
 \hoffset=1cm
 \voffset=2cm

\footline={\hss{\vbox to 2cm{\vfil\hbox{\rm\folio}}}\hss}
\nopagenumbers
\def\DJ{\leavevmode\setbox0=\hbox{D}\kern0pt\rlap
{\kern.04em\raise.188\ht0\hbox{-}}D}

\def\txt#1{{\textstyle{#1}}}
\baselineskip=13pt
\def\hf{{\textstyle{1\over2}}}
\def\a{\alpha}
\def\d{{\,\roman d}}
\def\e{\varepsilon}
\def\f{\varphi}
\def\G{\Gamma}

\def\s{\sigma}
\def\t{\theta}
\def\={\;=\;}

\def\zt{\zeta(\hf+it)}

\def\D{\Delta}
\def\no{\noindent}
 
\def\z{\zeta}

 \def\t{\theta}
\def\hf{{\textstyle{1\over2}}}
\def\txt#1{{\textstyle{#1}}}
\def\f{\varphi}

\font\tenmsb=msbm10
\font\sevenmsb=msbm7
\font\fivemsb=msbm5
\newfam\msbfam
\textfont\msbfam=\tenmsb
\scriptfont\msbfam=\sevenmsb
\scriptscriptfont\msbfam=\fivemsb
\def\Bbb#1{{\fam\msbfam #1}}

\def \RR {\Bbb R}

\font\ff=cmr8
\def\txt#1{{\textstyle{#1}}}
\baselineskip=13pt

\font\teneufm=eufm10
\font\seveneufm=eufm7
\font\fiveeufm=eufm5
\newfam\eufmfam
\textfont\eufmfam=\teneufm
\scriptfont\eufmfam=\seveneufm
\scriptscriptfont\eufmfam=\fiveeufm
\def\mathfrak#1{{\fam\eufmfam\relax#1}}

\font\tenmsb=msbm10
\font\sevenmsb=msbm7
\font\fivemsb=msbm5
\newfam\msbfam
     \textfont\msbfam=\tenmsb
      \scriptfont\msbfam=\sevenmsb
      \scriptscriptfont\msbfam=\fivemsb
\def\Bbb#1{{\fam\msbfam #1}}

\def \RR {\Bbb R}

  \def\rightheadline{{\hfil{\ff
  On the Riemann zeta-function and the divisor problem IV}\hfil\tenrm\folio}}

  \def\leftheadline{{\tenrm\folio\hfil{\ff
   Aleksandar Ivi\'c }\hfil}}
  \def\emptyheadline{\hfil}
  \headline{\ifnum\pageno=1 \emptyheadline\else
  \ifodd\pageno \rightheadline \else \leftheadline\fi\fi}

\topmatter
\title
ON THE RIEMANN ZETA-FUNCTION AND THE DIVISOR PROBLEM IV
\endtitle
\author   Aleksandar Ivi\'c  \endauthor
\address
Aleksandar Ivi\'c, Katedra Matematike RGF-a
Universiteta u Beogradu, \DJ u\v sina 7, 11000 Beograd, Serbia
\endaddress
\email {\tt
ivic\@rgf.bg.ac.yu,  aivic\@matf.bg.ac.yu} \endemail
\dedicatory
Uniform Distribution Theory 1(2006), 125-135.
\enddedicatory
\abstract
{Let $\D(x)$ denote the error term in the Dirichlet
divisor problem, and $E(T)$ the error term in the asymptotic
formula for the mean square of $|\zt|$. If
$E^*(t) = E(t) - 2\pi\D^*(t/2\pi)$ with $\D^*(x) =
 -\D(x)  + 2\D(2x) - \hf\D(4x)$, then it is proved that
$$
\int_0^T|E^*(t)|^3\d t \ll_\e T^{3/2+\e},
$$
which is (up to `$\e$') best possible,
and $\zt \ll_\e t^{\rho/2+\e}$ if $E^*(t) \ll_\e t^{\rho+\e}$.
}
\endabstract
\endtopmatter

\document
\head
1. Introduction and statement of results
\endhead

This paper is the continuation of the author's works [5], [6], where the analogy
between the Riemann zeta-function $\z(s)$ and the divisor problem
was investigated. As usual, let the error term in the
classical Dirichlet divisor problem be
$$
\D(x) \;=\; \sum_{n\le x}d(n) - x(\log x + 2\gamma - 1),
\leqno(1.1)
$$
and the error term in the mean square formula for $|\zt|$ be defined by
$$
E(T) \;=\;\int_0^T|\zt|^2\d t - T\left(\log\bigl({T\over2\pi}\bigr) + 2\gamma - 1
\right).\leqno(1.2)
$$
Here, as usual, $d(n)$ is the number of divisors of
$n$, $\z(s)$ is the Riemann zeta-function, and $ \gamma = -\G'(1) = 0.577215\ldots\,$
is Euler's constant. The analogy between $\zeta(s)$ and the divisor problem
is more exact if, instead with $\D(x)$, we work with the modified
function $\D^*(x)$ (see  M. Jutila [8], [9] and T. Meurman [11], [12]), where
$$
\D^*(x) := -\D(x)  + 2\D(2x) - \hf\D(4x)
= \hf\sum_{n\le4x}(-1)^nd(n) - x(\log x + 2\gamma - 1).
\leqno(1.3)
$$
M. Jutila (op. cit.) investigated both the
local and global behaviour of the difference
$$
E^*(t) \;:=\; E(t) - 2\pi\D^*\bigl({t\over2\pi}\bigr).\leqno(1.4)
$$
This function may be thought of as a discrepancy between $E^*(t)$ and $\D^*(x)$.
In particular Jutila in [9]  proved that
$$
\int_0^T(E^*(t))^2\d t \;\ll\; T^{4/3}\log^3T,\leqno(1.5)
$$
which was sharpened in [6] by the author  to the full asymptotic formula
$$
\int_0^T (E^*(t))^2\d t \;=\; T^{4/3}P_3(\log T) + O_\e(T^{7/6+\e}),\leqno(1.6)
$$
where $P_3(y)$ is a polynomial of degree three in $y$ with
positive leading coefficient, and all the coefficients may be evaluated
explicitly.
Here and later $\e$ denotes positive constants which are arbitrarily
small, but are not necessarily the same ones at each occurrence,
while $a \ll_\e b$ (same as $a = O_\e(b))$ means that
the $\ll$--constant depends on $\e$.
In Part II of [5] it was proved that
$$
\int_0^T |E^*(t)|^5\d t \;\ll_\e\; T^{2+\e},\leqno(1.7)
$$
while in Part III we investigated the function $R(T)$ defined by the relation
$$
\int_0^T E^*(t)\d t = {3\pi\over 4}T + R(T),\leqno(1.8)
$$
and proved, among other things, the asymptotic formula
$$
\int_0^TR^2(t)\d t = T^2Q_3(\log T) + O_\e(T^{11/6+\e}),\leqno(1.9)
$$
where $Q_3(y)$ is a cubic polynomial in $y$ with positive leading
coefficient, whose all coefficients may be  evaluated explicitly.

\medskip
The asymptotic formula (1.9) bears resemblance to (1.6), and it is proved
by a similar technique. The exponents in the error terms are, in both cases,
less than the exponent of $T$ in the main term by 1/6. This comes from the
use of [6, Lemma 3], and in both cases the exponent of the error
term is the limit of the method. Our first new result is an upper bound for the third moment
of $|E^*(t)|$, which does not follow from any of the previous results. This is

\bigskip
THEOREM 1. {\it We have}
$$
\int_0^T |E^*(t)|^3\d t \;\ll_\e\; T^{3/2+\e}.\leqno(1.10)
$$

\bigskip
In view of (1.6) it follows that, up to `$\e$', (1.10) is best possible.

\medskip
{\bf Corollary 1}. We have
$$
\int_0^T|\zt|^8\d t \;\ll_\e\; T^{3/2+\e}.
$$
\bigskip
The last result is, up to `$\e$', the sharpest one known
(see [3, Chapter 8]). It follows
from Theorem 1.4 of [5, Part II], which says that the bound
$$
\int_0^T|E^*(t)|^{k}\d t \;\ll_\e\; T^{c(k)+\e}\leqno(1.11)
$$
implies that
$$
\int_0^T|\zt|^{2k+2}\d t \;\ll_\e\; T^{c(k)+\e},\leqno(1.12)
$$
where $k\ge1$ is a fixed real number.

\medskip
{\bf Corollary 2}. We have
$$
\int_0^T(E^*(t))^4\d t \ll_\e T^{7/4+\e},\quad
\int_0^T|\zt|^{10}\d t \ll_\e T^{7/4+\e}.\leqno(1.13)
$$
\medskip
The first bound in (1.13) follows by the Cauchy-Schwarz inequality
for integrals from (1.7) and (1.10). The second bound follows from
(1.11)--(1.12) with $k=4$ and represents, up to `$\e$', the sharpest one known
(see [3, Chapter 8]). The first exponent in (1.13) improves on $16/9+\e$,
proved in [5, Part I].

\medskip
{\bf Corollary 3}. If, for $k>0$ a fixed constant and
$1 \ll G = G(T) \ll T$,
$$
J_k(T,G) := {1\over \sqrt{\pi}G}\int_{-\infty}^\infty
|\z(\hf + iT + iu)|^{2k}{\roman e}^{-(u/G)^2}\d u,
$$
then
$$
\int_T^{2T} J_1^4(t,G)\d t \;\ll_\e\;T^{1+\e} \leqno(1.14)
$$
holds for $T^{3/16} \le G = G(T) \ll T$.

\medskip
Namely it was proved in [6] that, for $\,T^\e \ll G = G(T) \le T$
and fixed $m\ge 1$ we have
$$
\int_T^{2T}J_1^m(t,G)\d t \ll G^{-1-m}\int_{-G\log T}^{G\log T}
\left(\int_T^{2T}|E^*(t+x)|^m\d t\right)\d x + T\log^{2m}T.\leqno(1.15)
$$
Thus (1.14) follows from (1.13) and (1.15) with $m=4$, and improves on the
range $T^{7/36} \le G = G(T) \ll T$ stated in Theorem 1 of [6], since $3/16
< 7/36$.

\medskip
Both (1.6) and (1.10) imply that, in the mean sense, $E^*(t) \ll_\e t^{1/6+\e}$.
The true order of this function is, however, quite elusive. If we define
$$
\rho \;:=\; \inf\Bigl\{\;r \,>0\;:\; E^*(T) = O(T^r)\;\Bigr\},\leqno(1.16)
$$
then we have unconditionally
$$
1/6 \;\le\; \rho\le 131/416 = 0.314903\ldots\,,\leqno(1.17)
$$
and there is a big discrepancy between the lower and upper bound in (1.17).
The lower bound in (1.17) comes from the asymptotic formula (1.6),
which in fact gives $E^*(T)
= \Omega(T^{1/6}(\log T)^{3/2})$. The
upper bound comes from the best known bound for $\D(x)$ of M.N. Huxley [2]
and $E(T)$ of N. Watt (unpublished). It remains yet to see whether a method
can be found that would provide sharper bounds for $\rho$ than for the
corresponding exponents of $E(T)$ and $\D(x)$. This is important, as one can
obtain bounds for $\zt$ from bounds of $E^*(t)$. More precisely, if as
usual one defines the Lindel\"of function for $\z(s)$ (the famous Lindel\"of conjecture
is that $\mu(\hf) = 0$) by the relation
$$
\mu(\s) \;=\; \liminf_{t\to\infty}\,{\log|\z(\s+it)|\over\log t}\leqno(1.18)
$$
for any $\s\in\RR$, then we have

\bigskip
THEOREM 2. {\it If $\rho$ is defined by} (1.16) {\it and $\mu(\s)$ by} (1.18),
{\it then we have}
$$
\mu(\hf) \;\le\; \hf\rho.\leqno(1.19)
$$

\bigskip
It may be remarked that, if $\rho\le 1/4$ holds, then
$\t = \omega$, where
$$
\t = \inf\Bigl\{\; c> 0\;:\;E(T) = O(T^c)\;\Bigr\},\quad
\omega = \inf\Bigl\{\; d > 0\;:\;\D(T) = O(T^d)\;\Bigr\}.
$$
Namely as $\t\ge1/4$ and $\omega\ge1/4$ are known to hold (this follows e.g., from
mean square results, see [4]) $\t = \omega$ follows from (1.4) and
$\omega = \sigma$, proved recently by Lau--Tsang [10], where
$$
\s = \inf\Bigl\{\; s> 0\;:\;\D^*(T) = O(T^s)\;\Bigr\}.
$$
The reader is also referred to M. Jutila [8] for a discussion on some related
implications. The limit of (1.19) is $\mu(\hf) \le 1/12$ in view of (1.17).
\bigskip
The plan of the paper is as follows. In Section 2 the necessary lemmas are
given, while the proofs of Theorem 1 and Theorem 2 will be given in
Section 3.

\smallskip

\head
2. The necessary lemmas
\endhead

\bigskip
In this section we shall state the lemmas which are necessary
for the proof of our theorems.

\medskip
LEMMA 1 (O. Robert--P. Sargos [13]). {\it Let $k\ge 2$ be a fixed
integer and $\delta > 0$ be given.
Then the number of integers $n_1,n_2,n_3,n_4$ such that
$N < n_1,n_2,n_3,n_4 \le 2N$ and}
$$
|n_1^{1/k} + n_2^{1/k} - n_3^{1/k} - n_4^{1/k}| < \delta N^{1/k}
$$
{\it is, for any given $\e>0$,}
$$
\ll_\e N^\e(N^4\delta + N^2).\leqno(2.1)
$$

This Lemma (with $k=2$) is crucial in treating the
fourth power of the sums in (2.5) and (2.12).
\medskip
LEMMA 2. {\it Let $T^\e \ll G \ll T/\log T$. Then we have}
$$
E^*(T) \le {2\over\sqrt{\pi}G}\int_0^{\infty} E^*(T+u)\,
{\roman e}^{-u^2/G^2}\d u + O_\e(GT^\e),\leqno(2.2)
$$
{\it and}
$$
E^*(T) \ge {2\over\sqrt{\pi}G}\int_0^\infty E^*(T-u)\,{\roman e}^{-u^2/G^2}\d u
+ O_\e(GT^\e).\leqno(2.3)
$$

\medskip \no
Lemma 2 follows on combining Lemma 2.2 and Lemma 2.3 of [4, Part I].

\medskip
The next lemma is F.V. Atkinson's classical, precise asymptotic
formula for $E(T)$ (see [1], [3] or [4]).

\medskip
LEMMA 3. {\it Let $0 < A < A'$ be any two fixed constants
such that $AT < N < A'T$, and let $N' = N'(T) =
T/(2\pi) + N/2 - (N^2/4+ NT/(2\pi))^{1/2}$. Then }
$$
E(T) = \Sigma_1(T) + \Sigma_2(T) + O(\log^2T),\leqno(2.4)
$$
{\it where}
$$
\Sigma_1(T) = 2^{1/2}(T/(2\pi))^{1/4}\sum_{n\le N}(-1)^nd(n)n^{-3/4}
e(T,n)\cos(f(T,n)),\leqno(2.5)
$$
$$
\Sigma_2(T) = -2\sum_{n\le N'}d(n)n^{-1/2}(\log T/(2\pi n))^{-1}
\cos(T\log T/(2\pi n) - T + \pi /4),\leqno(2.6)
$$
{\it with}
$$
\eqalign{\cr&
f(T,n) = 2T{\roman {arsinh}}\,\bigl(\sqrt{\pi n/(2T})\bigr) + \sqrt{2\pi nT
+ \pi^2n^2} - \pi/4\cr&
=  -\txt{1\over4}\pi + 2\sqrt{2\pi nT} +
\txt{1\over6}\sqrt{2\pi^3}n^{3/2}T^{-1/2} + a_5n^{5/2}T^{-3/2} +
a_7n^{7/2}T^{-5/2} + \ldots\,,\cr}\leqno(2.7)
$$
$$\eqalign{\cr
e(T,n) &= (1+\pi n/(2T))^{-1/4}{\Bigl\{(2T/\pi n)^{1/2}
{\roman {arsinh}}\,(\sqrt{\pi n/(2T})\Bigr\}}^{-1}\cr&
= 1 + O(n/T)\qquad(1 \le n < T),
\cr}\leqno(2.8)
$$
{\it and $\,{\roman{arsinh}}\,x = \log(x + \sqrt{1+x^2}\,).$}

\medskip
LEMMA 4 (M. Jutila [8, Part II]). {\it For $A\in\RR$ a constant we have}
$$
\cos\left(\sqrt{8\pi nT} +
\txt{1\over6}\sqrt{2\pi^3}n^{3/2}T^{-1/2} + A\right)
= \int_{-\infty}^\infty \a(u)\cos(\sqrt{8\pi n}(\sqrt{T} + u)
+ A)\d u,\leqno(2.9)
$$
{\it where $\a(u) \ll T^{1/6}$ for $u\not=0$,}
$$
\a(u) \ll T^{1/6}\exp(-bT^{1/4}|u|^{3/2}) \leqno(2.10)
$$
{\it for $u<0$, and}
$$
\a(u) = T^{1/8}u^{-1/4}\left(d\exp(ibT^{1/4}u^{3/2})
+ {\bar d}\exp(-ibT^{1/4}u^{3/2})\right) + O(T^{-1/8}u^{-7/4})\leqno(2.11)
$$
{\it for $u \ge T^{-1/6}$ and some constants $b\; (>0)$ and $d$.}

\bigskip
We need also an explicit formula for $\D^*(x)$ (see [3, Chapter 15]). This is

\bigskip
LEMMA 5. {\it For $1 \le N \ll x$ we have}
$$
\D^*(x) = {1\over\pi\sqrt{2}}x^{1\over4}
\sum_{n\le N}(-1)^nd(n)n^{-{3\over4}}
\cos(4\pi\sqrt{nx} - {\txt{1\over4}}\pi) +
O_\e(x^{{1\over2}+\e}N^{-{1\over2}}).
\leqno(2.12)
$$

\bigskip

\head
3. Proofs of the theorems
\endhead
\medskip
The proof of (1.10) of Theorem 1 is based on the method of [5].
We seek an upper bound for $R = R(V,T)$, the number of  points
$$
\{t_r\}\in[T,2T]\,(r = 1,\ldots,\, R),
\quad V \le |E^*(t_r)| < 2V \quad(|t_r-t_s|\ge V \;\roman {if}\;
r\ne s).\leqno(3.1)
$$
We consider separately the points where $E^*(t_r)$ is
positive or negative. Suppose the first case holds (the other one
is treated analogously), using in either case the notation $R$ for
the number of points in question. Then from Lemma 2 we have
$$ V \le E^*(t_r) \le
{2\over\sqrt{\pi}G}\int_0^\infty  E^*(t_r+G+u) \,{\roman
e}^{-u^2/G^2}\d u + O_\e(GT^\e),\leqno(3.2)
$$
and the integral may be truncated at $u = G\log T$ with a very small error.
We may suppose that $V$ satisfies
$$
T^{1/6} \;\le \;V \;\le T^{1/4}.\leqno(3.3)
$$
Indeed, if
$$
I_1(T) := \int_{T,|E^*|\le T^{1/6}}^{2T}|E^*(t)|^3\d t,\quad
I_2(T) := \int_{T,|E^*|\ge T^{1/4}}^{2T}|E^*(t)|^3\d t,
$$
then from (1.6) it follows that
$$
I_1(T) \le T^{1/6}\int_T^{2T}|E^*(t)|^2\d t \ll T^{3/2}\log^3T,\leqno(3.4)
$$
while from (1.7) we obtain that
$$
I_2(T) \le T^{-1/2}\int_T^{2T}|E^*(t)|^5\d t \ll_\e T^{3/2+\e}.\leqno(3.5)
$$
Thus supposing that (3.3) holds we estimate
$$
I(V,T) := \int_{T,V\le|E^*(t)|\le2V}^T|E^*(t)|^3\d t
$$
by splitting the interval $[T,\,2T]\,$ into $R\,(= R(V,T))$ disjoint
subintervals $J_r$ of length $\le V$, where in the $r$-th of these
intervals we define $t_r\,(r =1,\ldots\,,R)$ by
$$
|E^*(t_r)| = \sup_{t\in J_r}|E^*(t)|.
$$
The proof of Theorem 1 will be a consequence of the bound
$$
R \;\ll_\e\;T^{3/2+\e}V^{-4},\leqno(3.6)
$$
provided that (3.1) holds (considering separately points with even and odd
indices so that $|t_r-t_s|\ge V\,(r\ne s)$ is satisfied). Namely we have
$$
I(V,T) \ll V\sum_{j=1}^R|E^*(t_r)|^3 \ll_\e VT^{3/2+\e}V^{-4}V^3
 = T^{3/2+\e},\leqno(3.7)
$$
and from (3.4), (3.5) and (3.7) we obtain
$$
\int_{T}^{2T}|E^*(t)|^3\d t \;\ll_\e\; T^{3/2+\e}.\leqno(3.8)
$$
The bound (1.10) follows from (3.8) if one replaces $T$ by $T2^{-j}$ and
sums the corresponding results for $j = 1,2,\ldots\,$.

\medskip
We continue the proof of Theorem 1 by noting that,
like in [5, Part I], the integral on the
right-hand side of (3.2) is simplified by Atkinson's
formula (Lemma 3) and the truncated formula for $\D^*(x)$ (Lemma 5).
We take $G = cVT^{-\e}$ (with sufficiently small $c>0$)
to make the $O$-term in (3.2) $\le \hf V$, and then we obtain
$$
V \,\ll\, \sum_{j=4}^6 V^{-1}T^\e\int_0^{G\log T}\sum\nolimits_j(t_r+G+u)\,{\roman
e}^{-u^2/G^2}\d u\quad(r = 1, \ldots,\, R),  \leqno(3.9)
$$
where we choose $X = T^{1/3-\e},\, N = TG^{-2}\log T$ and, similarly to [5],
for $t\asymp T$ we set (in the notation of Lemma 3)
$$\eqalign{\cr&
\sum\nolimits_4(t) := t^{1/4}\sum_{X<n\le N}
(-1)^n d(n)n^{-3/4}e(t+u,n)\cos(f(t+u,n)),\cr&
\sum\nolimits_5(t) := t^{1/4}
\sum_{X<n\le N}(-1)^n d(n)n^{-3/4}\cos(\sqrt{8\pi n(t+u)}-\pi/4),
\cr}\leqno(3.10)
$$
$$
\sum\nolimits_6(t) := t^{-1/4}\sum_{n\le X}
(-1)^n d(n)n^{3/4} \cos(\sqrt{8\pi n(t+u)}-\pi/4).\leqno(3.11)
$$
The sums in (3.10)--(3.11) over $n$ are split into $O(\log T)$ subsums
over the ranges $K<n\le K'\le 2K$. We denote these sums by $\Sigma_j(t,K)$
and let
$\f(t)$ denote a smooth, nonnegative function supported in
$\,[T/2,\,5T/2]\,$, such that $\f(t) = 1$ when $T \le t \le 2T$.
There must exist a set of $M = M(K)$ points $\{\tau_m\}\in \{t_r\}$
such that $M(K) \gg R/\log T$ for some $j,K$, so that it suffices to
majorize $M(K)$, which we shall (with a slight abuse of notation)
henceforth denote again by $R$.
The contribution of $\sum_6(t,K)$ is estimated by raising
the relevant portion of (3.9) to the fourth power and summing over $r$,
noting that $|t_r-t_s|\ge V\,(r\ne s)$, so that the sum of integrals
over the intervals $[t_r +G, t_r + G + G\log T]$
is majorized by the integral over $[T/2,\,5T/2]$.
We proceed as in [5, Part I and Part II] integrating by parts, and using
$\f^{(\ell)}(t) \ll_\ell T^{-\ell}\;(\ell \ge0)$. It transpires, when we develop
$\sum^4_6(t,K)$ and set
$$
\D := \sqrt{n_1} + \sqrt{n_2} - \sqrt{n_3} - \sqrt{n_4}\,,
$$
that the contribution of $\D \ge T^{\e-1/2}$ is negligible (i.e., it is
smaller than
$T^{-A}$ for any given $A>0$). The contribution of $\D < T^{\e-1/2}$ is treated
by Lemma 1 and trivial estimation of the ensuing integral. We obtain
$$\eqalign{
 RV^4 &\ll V^{-1}T^\e\sup_{|u|\le G\log T}
\int_{T/2}^{2T} \f(t)\sum\nolimits^4_6(t,K)\,\d t\cr&
\ll_\e T^{1+\e}V^{-1}\sup_{|u|\le G\log T,|\D|\le T^{\e-1/2}}T^{-1}K^3(K^4K^{-1/2}|\D|
+ K^2)\cr&
\ll_\e T^\e V^{-1}(T^{-1/2}X^{13/2} + X^5) \ll_\e T^{5/3+\e}V^{-1},\cr}
$$
since $K \ll X = T^{1/3-\e}$. This gives, since (3.3) holds,
$$
R \ll_\e T^{5/3+\e}V^{-5} \ll_\e T^{3/2+\e}V^{-4},
$$
which is the desired bound (3.6).

\medskip
The contributions of $\sum_4(t,K)$ and of $\sum_5(t,K)$
are estimated analogously, with the remark that in the case of
$\sum_4(t,K)$ one has to use Lemma 4 to deal
with the complications arising from the presence of $\cos(f(t+u,n))$,
coming from (2.5). This procedure was explained in detail
in [5, Part I and Part II]. The non-negligible
contribution of $\sum_5(t,K)$ will, again
by raising the relevant expression to the fourth power,
be  for $\D \le T^{\e-1/2}$ again. The application
of Lemma 1  gives in this case
$$
\eqalign{
RV^4 &\ll_\e V^{-1}T^{1+\e}TK^{-3}(K^4T^{-1/2} + K^2)\cr&
\ll_\e T^{2+\e}V^{-1}(K^{1/2}T^{1/2} + K^{-1})\cr&
\ll_\e T^{3/2+\e}V^{-1}K^{1/2} + T^{5/3+\e}V^{-1},\cr}\leqno(3.12)
$$
because $K\gg X = T^{1/3-\e}$ holds. For $K\le V^2$ the bound (3.12) reduces
to (3.6), and we are done. If $V^2 < K \le T^{1+\e}V^{-2}$ (note that
$V^2 <  T^{1+\e}V^{-2}$ holds by (3.3)), then the relevant expression
is squared, and not raised to the fourth power. We obtain
$$
\eqalign{
RV^2 &\ll_\e V^{-1}\max_{|u|\le G\log T}\int_{T/2}^{5T/2}\f(t)
\sum\nolimits_5^2(t,K)\d t\cr&
= T^{1/2}V^{-1}\max_{|u|\le G\log T}\int_{T/2}^{5T/2}\f(t)\times\cr&
\times
\sum_{K<m,n\le2K}(-1)^{m+n}d(m)d(n)(mn)^{-3/4}
{\roman e}^{i\sqrt{8\pi(t+u)}(\sqrt{m}-\sqrt{n})}\d t\cr&
\ll T^{3/2}V^{-1}\sum_{m>K}d^2(m)m^{-3/2} +
T^{1+\e}K^{-3/2}V^{-1}\sum_{K<m\ne n\le 2K}|\sqrt{m}-\sqrt{n}|^{-1}.
\cr}
$$
Here we used trivial estimation for the diagonal terms $m=n$, and the first
derivative test ([3, Lemma 2.1]) for the remaining terms. Since $V^2 < K$
and
$$
\sum_{K<m\ne n\le 2K}|\sqrt{m}-\sqrt{n}|^{-1}
\ll \sum_{K<m\le2K} \sqrt{K}\sum_{K<n\le2K, n\ne m}|m-n|^{-1}
\ll K^{3/2}\log K,
$$
we obtain that
$$
RV^2 \ll_\e T^{3/2}V^{-1}K^{-1/2}\log^3T + T^{1+\e}V^{-1} \ll_\e T^{3/2+\e}V^{-2},
$$
and (3.6) follows again. The proof of Theorem 1 is complete.

\bigskip
For the proof of Theorem 2 note that,
by [4, Theorem 1.2], (1.4)  and (1.19), we have
$$
\eqalign{&
|\z(\hf + iT)|^2 \ll \log T \int_{T-1}^{T+1}|\zt|^2\d t + 1\cr&
\ll \log T\Bigl(\log T + E(T+1) - E(T-1)\Bigr)\cr&
\ll_\e \log T\left(\log T + 2\pi\D^*\Bigl({T+1\over2\pi}\Bigr) -
2\pi\D^*\Bigl({T-1\over2\pi}\Bigr)\right) + T^{\rho+\e} \ll_\e T^{\rho+\e},\cr}
\leqno(3.13)
$$
since, from (1.3) and $d(n) \ll_\e n^\e$, it is seen that
$$
\D^*(T+H) - \D^*(T) = O(H\log T) +\hf\sum_{4T<n\le4(T+H)}(-1)^nd(n) \ll_\e HT^\e
$$
holds for $1 \ll H \ll T$. Therefore (3.13) implies that
$$
|\z(\hf + iT)|^2 \;\ll_\e\; T^{\rho+\e},
$$
and this gives $\mu(\hf) \le \hf\rho$, as asserted.

\bigskip
\Refs
\bigskip

\item{[1]} F.V. Atkinson, The mean value of the Riemann zeta-function,
Acta Math. {\bf81}(1949), 353-376.

\item{[2]} M.N. Huxley, Exponential sums and the Riemann zeta-function IV,
Proc. London Math. Soc. (3) {\bf66}(1993), 1-40 and V, ibid.
(3) {\bf90}(2005), 1-41.

\item{[3]} A. Ivi\'c, The Riemann zeta-function, John Wiley \&
Sons, New York, 1985 (2nd ed. Dover, Mineola, New York, 2003).

\item{[4]} A. Ivi\'c, The mean values of the Riemann zeta-function,
LNs {\bf 82}, Tata Inst. of Fundamental Research, Bombay (distr. by
Springer Verlag, Berlin etc.), 1991.

\item{[5]} A. Ivi\'c, On the Riemann zeta-function and the divisor problem,
Central European J. Math. {\bf(2)(4)} (2004), 1-15; II, ibid.
{\bf(3)(2)} (2005), 203-214, and III, subm. to Ann. Univ. Budapest. Sectio
Computatorica, also {\tt arXiv:math.NT/0610539}.

\item{[6]} A. Ivi\'c, On the mean square of the zeta-function and
the divisor problem, Annales  Acad. Scien. Fennicae Mathematica (in press),
also {\tt arXiv:math.NT/0603491}.

\item{[7]} A. Ivi\'c, Some remarks on the moments of $|\zt|$ in short
intervals, subm. to Acta Math. Hungarica, also {\tt arXiv:math.NT/0611427}.

\item{[8]} M. Jutila, Riemann's zeta-function and the divisor problem,
Arkiv Mat. {\bf21}(1983), 75-96 and II, ibid. {\bf31}(1993), 61-70.

\item{[9]} M. Jutila, On a formula of Atkinson, in ``Coll. Math. Sci.
J\'anos Bolyai 34, Topics in classical Number Theory, Budapest 1981",
North-Holland, Amsterdam, 1984, pp. 807-823.

\item{[10]} Y.-K. Lau and K.-M. Tsang, Omega result for the mean square of the
Riemann zeta-function, Manuscripta Math. {\bf117}(2005), 373-381.

\item{[11]} T. Meurman, A generalization of Atkinson's formula to $L$-functions,
Acta Arith. {\bf47}(1986), 351-370.

\item{[12]} T. Meurman, On the mean square of the Riemann zeta-function,
Quart. J. Math. (Oxford) {\bf(2)38}(1987), 337-343.

\item{[13]} O. Robert and P. Sargos, Three-dimensional
exponential sums with monomials, J. reine angew. Math. {\bf591}(2006), 1-20.

\endRefs
\vskip1cm

\enddocument

\end